\documentclass{amsart}%
\usepackage{amsfonts}
\usepackage{amsmath}
\usepackage{amssymb}
\usepackage{graphicx}
\usepackage{subfig}
\usepackage{eurosym}
\usepackage{amssymb}
\usepackage{amsmath}
\usepackage{amsfonts}
\usepackage{graphicx}%
\providecommand{\U}[1]{\protect\rule{.1in}{.1in}}
\newtheorem{theorem}{Theorem}
\theoremstyle{plain}

\newtheorem{definition}{Definition}
\newtheorem{example}{Example}

\numberwithin{equation}{section}
\begin{document}
\title[Firat University]{Assesment of Smarandache Curves in The Null Cone $\mathbf{Q}^{2}$}
\author{Mihriban Kulahci}
\address{Department of Mathematics, Firat University, 23119 ELAZI\u{G}/T\"{U}RK\.{I}YE}
\email{mihribankulahci@gmail.com}
\author{Fatma Almaz}
\email{fb\_fat\_almaz@hotmail.com}
\thanks{This paper is in final form and no version of it will be submitted for
publication elsewhere.}
\date{2017}

\begin{abstract}
In this paper, we give Smarandache curves according to the asymptotic
orthonormal frame in null cone $\mathbf{Q}^{2}$. By using cone frame formulas,
we present some characterizations of Smarandache curves and calculate cone
frenet invariants of these curves. Also, we illustrate these curves with an example.

\end{abstract}
\maketitle

\section{Introduction}

Human being were bewitched by curves and curved shapes long before they took
into account them as mathematical objects. But the greatest effect in the
research of curves was, of course, the discovery of the calculus. Geometry
before calculus includes only the simplest curves.

In classical curve theory, the geometry of a curve in three-dimensions is
actually characterized by Frenet vectors.

Smarandache geometry is a geometry which has at least one Smarandachely denied
axiom \cite{4}. An axiom is said to be Smarandachely denied, if it behaves in
at least two different ways within the same space. Smarandache curve is
defined as a regular curve whose position vector is composed by Frenet frame
vectors of another regular curve. Smarandache curves in various ambient spaces
have been classfied in \cite{1}-\cite{14}, \cite{19}, \cite{21}-\cite{32}.

In this study, we define special Smarandache curves such as $x\alpha,xy,\alpha
y$ and $x\alpha y$ Smarandache curves according to asymptotic orthonormal
frame in the null cone $\mathbf{Q}^{2}$ and we examine the curvature and the
asymptotic orthonormal frame's vectors of Smarandache curves. We also give an
example related to these curves.

\section{\textbf{Preliminaries}}

Some basics of the curves in the null cone are provided from, \cite{15}%
-\cite{18}.

Let $E_{1}^{3}$ be the $3-$dimensional pseudo-Euclidean space with the%

\[
\widetilde{g}(X,Y)=\langle X,Y\rangle=x_{1}y_{1}+x_{2}y_{2}-x_{3}y_{3}%
\]
for all $X=(x_{1},x_{2},x_{3}),$ $Y=(y_{1},y_{2},y_{3})\in E_{1}^{3}$.
$E_{1}^{3}$ is a flat pseudo-Riemannian manifold of signature $(2,1)$.

Let $M$ be a submanifold of $E_{1}^{3}$. If the pseudo-Riemannian metric
$\widetilde{g}$ of $E_{1}^{3}$ induces a pseudo-Riemannian metric
$g$(respectively, a Riemannian metric, a degenerate quadratic form) on $M$,
then $M$ is called a timelike( respectively, spacelike, degenerate)
submanifold of $E_{1}^{3}.$ Let $c$ be a fixed point in $E_{1}^{3}.$ The
pseudo-Riemannian lightlike cone(quadric cone ) is defined by%

\[
\mathbf{Q}_{1}^{2}(c)=\left\{  x\in E_{1}^{3}:g(x-c,x-c)=0\right\}  ,
\]
where the point $c$ is called the center of $\mathbf{Q}_{1}^{2}(c)$. When
$c=0$, we simply denote $\mathbf{Q}_{1}^{2}(0)$ by $\mathbf{Q}_{{}}^{2}$ be
and call it the null cone.

Let $E_{1}^{3}$ be $3$-dimensional Minkowski space and $\mathbf{Q}^{2}$ the
lightlike cone in $E_{1}^{3}.$ A vector $V\neq0$ $\ $in $E_{1}^{3}$ is called
spacelike, timelike or lightlike, if $\langle V,V\rangle>0$, $\langle
V,V\rangle<0$ or $\langle V,V\rangle=0,$ respectively. The norm of a vector
$x\in E_{1}^{3}$ is given by $\left\Vert x\right\Vert =\sqrt{\left\langle
x,x\right\rangle },$ \cite{20}.

We assume that curve $x:I\rightarrow\mathbf{Q}_{{}}^{2}\subset E_{1}^{3}$ is a
regular curve in $\mathbf{Q}_{{}}^{2}$ for $t\in I.$ In the following, we
always assume that the curve is regular.

A frame field $\left\{  x,\alpha,y\right\}  $ on $E_{1}^{3}$ is called an
asymptotic orthonormal frame field, if%
\[
\langle x,x\rangle=\langle y,y\rangle=\langle x,\alpha\rangle=\langle
y,\alpha\rangle=0,\text{ }\langle x,y\rangle=\langle\alpha,\alpha\rangle=1.
\]

Using \ $x^{\prime}(s)=\alpha(s)$ we know that $\left\{  x(s),\alpha
(s),y(s)\right\}  $ from an asymptotic orthonormal frame along the curve
$x(s)$ and the cone frenet formulas of $x(s)$ are given by%

\begin{align}
x^{\prime}(s)  &  =\alpha(s)\nonumber\\
\alpha^{\prime}(s)  &  =\kappa(s)x(s)-y(s)\tag{2.1}\label{2.1}\\
y^{\prime}(s)  &  =-\kappa(s)\alpha(s)\nonumber
\end{align}
where the function $\kappa(s)$ is called cone curvature function of the curve
$x(s)$, \cite{17}$.$

Let $x:I\rightarrow\mathbf{Q}_{{}}^{2}\subset E_{1}^{3}$ be a spacelike curve
in $\mathbf{Q}_{{}}^{2}$ with an arc length parameter $s.$ Then $x=x(s)=(x_{1}%
,x_{2},x_{3})$ can be written as%

\begin{equation}
x(s)=\frac{f_{s}^{-1}}{2}(f^{2}-1,2f,f^{2}+1) \tag{2.2}\label{2.2}%
\end{equation}
for some non constant function $f(s)$ and $f_{s}=f^{\prime}$, \cite{18}$.$

\section{Smarandache Curves in The Null Cone $\mathbf{Q}^{2}$}

In this section, we define the Smarandache curves according to the asymptotic
orthonormal frame in $\mathbf{Q}_{{}}^{2}$. Also, we obtain the asymptotic
orthonormal frame and cone curvature function of the Smarandache partners
lying on $\mathbf{Q}_{{}}^{2}$ using cone frenet formulas.

Smarandache curve $\gamma=\gamma(s^{\ast}(s))$ of the curve $x$ is a regular
unit speed curve lying fully on $\mathbf{Q}_{{}}^{2}$. Let $\left\{
x,\alpha,y\right\}  $ and $\left\{  \gamma,\alpha_{\gamma},y_{\gamma}\right\}
$ be the moving asymptotic orthonormal frames of $x$ and $\gamma,$ respectively.

\begin{definition}
Let $x$ be unit speed spacelike curve lying on $\mathbf{Q}_{{}}^{2}$ with the
moving asymptotic orthonormal frame $\left\{  x,\alpha,y\right\}  .$ Then,
$x\alpha-$smarandache curve of $x$ is defined by%
\begin{equation}
\gamma_{x\alpha}(s^{\ast})=\frac{c}{b}x(s)+\alpha(s), \tag{3.1}\label{3.1}%
\end{equation}
where $c,b\in%
\mathbb{R}
_{0}^{+}.$

\begin{theorem}
Let $x$ be unit speed spacelike curve in $\mathbf{Q}_{{}}^{2}$ with the moving
asymptotic orthonormal frame $\left\{  x,\alpha,y\right\}  $ and cone
curvature $\kappa(s)$ and let $\gamma_{x\alpha}$ be $x\alpha-$smarandache
curve with asymptotic orthonormal frame $\left\{  \gamma_{x\alpha}%
,\alpha_{x\alpha},y_{x\alpha}\right\}  .$ Then the following relations hold:

i) The asymptotic orthonormal frame $\left\{  \gamma_{x\alpha},\alpha
_{x\alpha},y_{x\alpha}\right\}  $ of the $x\alpha$-smarandache curve
$\gamma_{x\alpha}$ is given as%
\begin{equation}%
\begin{bmatrix}
\gamma_{x\alpha}\\
\alpha_{x\alpha}\\
y_{x\alpha}%
\end{bmatrix}
=%
\begin{bmatrix}
\frac{c}{b} & 1 & 0\\
\psi & \frac{b}{c}\kappa\psi & -\frac{b}{c}\psi\\
\varrho_{1} & \varrho_{2} & \varrho_{3}%
\end{bmatrix}%
\begin{bmatrix}
x\\
\alpha\\
y
\end{bmatrix}
\tag{3.2}\label{3.2}%
\end{equation}
where
\begin{align}
\psi &  =\frac{c}{\sqrt{c^{2}-2\kappa b^{2}}},\tag{3.3}\label{3.3}\\
\Upsilon_{1}  &  =\frac{b}{c}(\frac{b}{c}\kappa^{\prime}+\kappa)\psi
\psi^{\prime}+\frac{b^{2}}{c^{2}}\kappa\left(  \psi^{\prime}\right)
^{2},\nonumber\\
\Upsilon_{2}  &  =\frac{b}{c}\psi(\psi^{\prime}+2\frac{b}{c}\kappa
\psi),\tag{3.4}\label{3.4}\\
\Upsilon_{3}  &  =-\frac{b}{c}\psi(\frac{b}{c}\psi^{\prime}+\psi)\nonumber
\end{align}
and
\begin{align}
\varrho_{1}  &  =-(\Upsilon_{1}+\frac{c}{2b}\left(  2\Upsilon_{1}\Upsilon
_{3}+\Upsilon_{2}^{2}\right)  )=-\Upsilon_{1}+\frac{c}{b}\kappa_{\gamma
_{x\alpha}}(s^{\ast}),\nonumber\\
\varrho_{2}  &  =-(\Upsilon_{2}+\frac{1}{2}\left(  2\Upsilon_{1}\Upsilon
_{3}+\Upsilon_{2}^{2}\right)  )=-\Upsilon_{2}+\kappa_{\gamma_{x\alpha}%
}(s^{\ast}),\tag{3.5}\label{3.5}\\
\varrho_{3}  &  =-\Upsilon_{3}.\nonumber
\end{align}

ii) The cone curvature $\kappa_{\gamma_{x\alpha}}(s^{\ast})$ of the curve
$\gamma_{x\alpha}$ is given by%
\begin{equation}
\kappa_{\gamma_{x\alpha}}(s^{\ast})=-\frac{1}{2}\left(  2\Upsilon_{1}%
\Upsilon_{3}+\Upsilon_{2}^{2}\right)  . \tag{3.6}\label{3.6}%
\end{equation}
where
\[
s^{\ast}=\frac{1}{b}\int\sqrt{c^{2}-2b^{2}\kappa(s)}ds.
\]

\end{theorem}
\end{definition}

\begin{proof}
\textbf{i)} We assume that the curve $x$ is a unit speed spacelike \ curve
with the asymptotic orthonormal frame $\left\{  x,\alpha,y\right\}  $ and cone
curvature $\kappa$. Differentiating the equation
\eqref{3.1}
with respect to $s$ and considering
\eqref{2.1}%
, we have%
\begin{equation}
\gamma_{x\alpha}^{\prime}(s^{\ast})=A\left(  \alpha(s)+\frac{b}{c}\kappa
x(s)-\frac{b}{c}y(s)\right)  ,\tag{3.7}\label{3.7}%
\end{equation}
where
\begin{align}
\frac{ds^{\ast}}{ds} &  =\frac{1}{b}\sqrt{c^{2}-2b^{2}\kappa(s)}%
,\tag{3.8}\label{3.8}\\
\psi &  =\frac{c}{\sqrt{c^{2}-2\kappa(s)b^{2}}}.\tag{3.9}\label{3.9}%
\end{align}

It can be easily seen that the tangent vector $\gamma_{x\alpha}^{\prime
}(s^{\ast})=\alpha_{x\alpha}(s^{\ast})$ is a unit spacelike vector.

Differentiating
\eqref{3.7}
, we obtain equation as follows%
\begin{equation}
\gamma_{x\alpha}^{\prime\prime}(s^{\ast})=\Upsilon_{1}x(s)+\Upsilon_{2}%
\alpha(s)+\Upsilon_{3}y(s), \tag{3.10}\label{3.10}%
\end{equation}
where%
\begin{align*}
\Upsilon_{1}  &  =\frac{b}{c}(\frac{b}{c}\kappa^{\prime}+\kappa)\psi
\psi^{\prime}+\frac{b^{2}}{c^{2}}\kappa\left(  \psi^{\prime}\right)  ^{2},\\
\Upsilon_{2}  &  =\frac{b}{c}\psi(\psi^{\prime}+2\frac{b}{c}\kappa\psi),\\
\Upsilon_{3}  &  =-\frac{b}{c}\psi(\frac{b}{c}\psi^{\prime}+\psi).
\end{align*}%
\begin{equation}
y_{x\alpha}(s^{\ast})=-\gamma_{x\alpha}^{\prime\prime}-\frac{1}{2}\left\langle
\gamma_{x\alpha}^{\prime\prime},\gamma_{x\alpha}^{\prime\prime}\right\rangle
\gamma_{x\alpha} \tag{3.11}\label{3.11}%
\end{equation}

By the help of previous equation
\eqref{3.11}%
, we obtain%
\begin{equation}
y_{x\alpha}(s^{\ast})=\varrho_{1}x(s)+\varrho_{2}\alpha(s)+\varrho_{3}y(s),
\tag{3.12}\label{3.12}%
\end{equation}
where $\varrho_{1}=-(\Upsilon_{1}+\frac{c}{2b}\left(  2\Upsilon_{1}%
\Upsilon_{3}+\Upsilon_{2}^{2}\right)  ),$ $\varrho_{2}=-(\Upsilon_{2}+\frac
{1}{2}\left(  2\Upsilon_{1}\Upsilon_{3}+\Upsilon_{2}^{2}\right)  ),$
$\varrho_{3}=-\Upsilon_{3}.$

\textbf{ii)} The curvature $\kappa_{\gamma_{x\alpha}}(s^{\ast})$ of the
$\gamma_{x\alpha}(s^{\ast})$ is explicity obtained by%
\begin{align}
\kappa_{\gamma_{x\alpha}}(s^{\ast})  &  =-\frac{1}{2}\left\langle
\gamma_{x\alpha}^{\prime\prime},\gamma_{x\alpha}^{\prime\prime}\right\rangle
\nonumber\\
&  =-\frac{1}{2}\left(  2\Upsilon_{1}\Upsilon_{3}+\Upsilon_{2}^{2}\right)  .
\tag{3.13}\label{3.13}%
\end{align}

Thus, the theorem is proved.
\end{proof}

\begin{definition}
Let $x$ be unit speed spacelike curve lying on $\mathbf{Q}_{{}}^{2}$ with the
moving asymptotic orthonormal frame $\left\{  x,\alpha,y\right\}  .$ Then,
$xy-$smarandache curve of $x$ is defined by%
\begin{equation}
\gamma_{xy}(s^{\ast})=\frac{1}{\sqrt{2cb}}\left(  cx(s)+by(s)\right)  ,
\tag{3.14}\label{3.14}%
\end{equation}
where $c,b\in%
\mathbb{R}
_{0}^{+}.$
\end{definition}

\begin{theorem}
Let $x$ be unit speed spacelike curve in $\mathbf{Q}_{{}}^{2}$ with the moving
asymptotic orthonormal frame $\left\{  x,\alpha,y\right\}  $ and cone
curvature $\kappa$ and let $\gamma_{xy}$ be $xy-$smarandache curve with
asymptotic orthonormal frame $\left\{  \gamma_{xy},\text{ }\alpha_{xy},\text{
}y_{xy}\right\}  .$ Then the following relations hold:

i) The asymptotic orthonormal frame $\left\{  \gamma_{xy},\alpha_{xy}%
,y_{xy}\right\}  $ of the $xy-$smarandache curve $\gamma_{xy}$ is given as%
\begin{equation}%
\begin{bmatrix}
\gamma_{xy}\\
\alpha_{xy}\\
y_{xy}%
\end{bmatrix}
=%
\begin{bmatrix}
\sqrt{\frac{c}{2b}} & 0 & \sqrt{\frac{b}{2c}}\\
0 & 1 & 0\\
\frac{b\kappa^{2}\sqrt{2bc}}{\left(  c-b\kappa\right)  ^{2}} & 0 &
\frac{c\sqrt{2bc}}{\left(  c-b\kappa\right)  ^{2}}%
\end{bmatrix}%
\begin{bmatrix}
x\\
\alpha\\
y
\end{bmatrix}
. \tag{3.15}\label{3.15}%
\end{equation}

ii) The cone curvature $\kappa_{\gamma_{xy}}(s^{\ast})$ of the curve
$\gamma_{xy}$ is given by%
\begin{equation}
\kappa_{\gamma_{xy}}(s^{\ast})=\frac{2bc\kappa(s)}{\left(  c-b\kappa
(s)\right)  ^{2}}, \tag{3.16}\label{3.16}%
\end{equation}
where
\begin{equation}
s^{\ast}=\frac{1}{\sqrt{2cb}}\int\left(  c-b\kappa(s)\right)  ds.
\tag{3.17}\label{3.17}%
\end{equation}

\end{theorem}

\begin{proof}
\textbf{i) }We assume that the curve $x$ is a unit speed spacelike \ curve
with the asymptotic orthonormal frame $\left\{  x,\alpha,y\right\}  $ and cone
curvature $\kappa$. Differentiating the equation
\eqref{3.14}
with respect to $s$ and considering
\eqref{2.1}%
, we have%
\begin{equation}
\gamma_{xy}^{\prime}(s^{\ast})\frac{ds^{\ast}}{ds}=\frac{1}{\sqrt{2cb}}\left(
c-b\kappa(s)\right)  \overrightarrow{\alpha}(s).\tag{3.18}\label{3.18}%
\end{equation}

By considering
\eqref{3.17}%
, we get%
\begin{equation}
\gamma_{xy}^{\prime}(s^{\ast})=\alpha(s)=\alpha_{xy}. \tag{3.19}\label{3.19}%
\end{equation}

Here, it can be easily seen that the tangent vector $\overrightarrow{\alpha
}_{xy}$ is a unit spacelike vector.%
\begin{equation}
\gamma_{xy}^{\prime\prime}(s^{\ast})\frac{ds^{\ast}}{ds}=\kappa x(s)-y(s).
\tag{3.20}\label{3.20}%
\end{equation}

By substituting
\eqref{3.17}
into
\eqref{3.20}
and making necessary calculations, we obtain%
\begin{equation}
\gamma_{xy}^{\prime\prime}(s^{\ast})=\frac{\kappa\sqrt{2bc}}{\left(
c-b\kappa\right)  }\overrightarrow{x}-\frac{\sqrt{2bc}}{\left(  c-b\kappa
\right)  ^{2}}\overrightarrow{y}. \tag{3.21}\label{3.21}%
\end{equation}

By the help of equation $y_{xy}(s^{\ast})=-\gamma_{xy}^{\prime\prime}-\frac
{1}{2}\left\langle \gamma_{xy}^{\prime\prime},\gamma_{xy}^{\prime\prime
}\right\rangle \gamma_{xy}$, we write%
\begin{equation}
y_{xy}(s^{\ast})=\frac{b\sqrt{2bc}.\kappa^{2}}{\left(  c-b\kappa\right)  ^{2}%
}x(s)+\frac{c\sqrt{2bc}}{\left(  c-b\kappa\right)  ^{2}}y(s). \tag{3.22}%
\label{3.22}%
\end{equation}

\textbf{ii)} The curvature $\kappa_{\gamma_{xy}}(s^{\ast})$ of the
$\gamma_{xy}(s^{\ast})$ is explicity obtained by%
\[
\kappa_{\gamma_{xy}}(s^{\ast})=-\frac{1}{2}\left\langle \gamma_{xy}%
^{\prime\prime},\gamma_{xy}^{\prime\prime}\right\rangle =\frac{2bc\kappa
(s)}{\left(  c-b\kappa(s)\right)  ^{2}}.
\]

\end{proof}

\begin{definition}
Let $x$ be unit speed spacelike curve lying on $\mathbf{Q}_{{}}^{2}$ with the
moving asymptotic orthonormal frame $\left\{  x,\alpha,y\right\}  .$ Then,
$\alpha y-$smarandache curve of $x$ is defined by%
\begin{equation}
\gamma_{\alpha y}(s^{\ast})=\alpha(s)+\frac{b}{c}y(s), \tag{3.23}\label{3.23}%
\end{equation}
where $c,b\in%
\mathbb{R}
_{0}^{+}.$
\end{definition}

\begin{theorem}
Let $x$ be unit speed spacelike curve in $\mathbf{Q}_{{}}^{2}$ with the moving
asymptotic orthonormal frame $\left\{  x,\alpha,y\right\}  $ and cone
curvature $\kappa$ and let $\gamma_{\alpha y}$ be $\alpha y-$smarandache curve
with asymptotic orthonormal frame $\left\{  \gamma_{\alpha y},\text{ }%
\alpha_{\alpha y},\text{ }y_{\alpha y}\right\}  .$ Then the following
relations hold:

i) The asymptotic orthonormal frame $\left\{  \gamma_{\alpha y},\alpha_{\alpha
y},y_{\alpha y}\right\}  $ of the $\alpha y-$smarandache curve $\gamma_{\alpha
y}$ is given as%
\begin{equation}%
\begin{bmatrix}
\gamma_{\alpha y}\\
\alpha_{\alpha y}\\
y_{\alpha y}%
\end{bmatrix}
=%
\begin{bmatrix}
0 & 1 & \frac{b}{c}\\
\frac{c\sqrt{\kappa}}{\sqrt{b^{2}-2c^{2}}} & \frac{b\sqrt{\kappa}}{\sqrt
{b^{2}-2c^{2}}} & \frac{c\sqrt{\kappa}}{\kappa\sqrt{b^{2}-2c^{2}}}\\
\omega_{1} & \omega_{2} & \omega_{3}%
\end{bmatrix}%
\begin{bmatrix}
x\\
\alpha\\
y
\end{bmatrix}
, \tag{3.24}\label{3.24}%
\end{equation}
where
\begin{align}
\zeta_{1}  &  =\frac{c\kappa^{\prime}}{b^{2}-2c^{2}}\left(  \frac{c-2\kappa
b}{2\kappa}\right)  ,\nonumber\\
\zeta_{2}  &  =\frac{c\kappa^{\prime}}{b^{2}-2c^{2}}\left(  \frac
{c-b-1}{2\kappa}\right)  ,\tag{3.25}\label{3.25}\\
\zeta_{3}  &  =\frac{c\kappa^{\prime}}{b^{2}-2c^{2}}\left(  \frac{b\kappa
-c}{2\kappa^{2}}\right) \nonumber
\end{align}
and
\begin{align}
\omega_{1}  &  =-\zeta_{1},\nonumber\\
\omega_{2}  &  =-(\zeta_{2}+\frac{1}{2}\left(  2\zeta_{1}\zeta_{3}+\zeta
_{2}^{2}\right)  ),\tag{3.26}\label{3.26}\\
\omega_{3}  &  =-(\zeta_{3}+\frac{b}{2c}\left(  2\zeta_{1}\zeta_{3}+\zeta
_{2}^{2}\right)  ).\nonumber
\end{align}

ii) The cone curvature $\kappa_{\gamma_{\alpha y}}(s^{\ast})$ of the curve
$\gamma_{\alpha y}$ is given by%
\begin{equation}
\kappa_{\gamma_{\alpha y}}(s^{\ast})=-\frac{c^{2}}{8\left(  b^{2}%
-2c^{2}\right)  ^{2}}(\frac{\kappa^{\prime}}{\kappa^{{}}})^{2}\left(  \frac
{1}{\kappa}\left(  c-2\kappa b\right)  \left(  b\kappa-c\right)  +\left(
c-b-1\right)  ^{2}\right)  , \tag{3.27}\label{3.27}%
\end{equation}
where
\begin{equation}
s^{\ast}=\frac{\sqrt{b^{2}-2c^{2}}}{c}\int\sqrt{\kappa(s)}ds. \tag{3.28}%
\label{3.28}%
\end{equation}

\end{theorem}

\begin{proof}
\textbf{i)} Let the curve $x$ be a unit speed spacelike\ curve with the
asymptotic orthonormal frame $\left\{  x,\alpha,y\right\}  $ and cone
curvature $\kappa$. Differentiating the equation
\eqref{3.23}
with respect to $s$ and considering
\eqref{2.1}%
, we find%
\[
\gamma_{\alpha y}^{\prime}(s^{\ast})\frac{ds^{\ast}}{ds}=\kappa\overrightarrow
{x(s)}-\frac{b}{c}\kappa\overrightarrow{\alpha(s)}-\overrightarrow{y(s)}.
\]

This can be written as following%
\begin{equation}
\alpha_{\alpha y}(s^{\ast})\frac{ds^{\ast}}{ds}=\kappa\overrightarrow
{x(s)}-\frac{b}{c}\kappa\overrightarrow{\alpha(s)}-\overrightarrow{y(s)},
\tag{3.29}\label{3.29}%
\end{equation}
where
\begin{equation}
\frac{ds^{\ast}}{ds}=\frac{\sqrt{b^{2}-2c^{2}}}{c}\sqrt{\kappa(s)}.
\tag{3.30}\label{3.30}%
\end{equation}

By substituting
\eqref{3.30}
into
\eqref{3.29}%
, we find%
\begin{equation}
\alpha_{\alpha y}(s^{\ast})=\frac{c\sqrt{\kappa}}{\sqrt{b^{2}-2c^{2}}%
}\overrightarrow{x}-\frac{b\sqrt{\kappa}}{\sqrt{b^{2}-2c^{2}}}\overrightarrow
{\alpha}-\frac{c\sqrt{\kappa}}{\kappa\sqrt{b^{2}-2c^{2}}}\overrightarrow{y}.
\tag{3.31}\label{3.31}%
\end{equation}

Differentiating
\eqref{3.31}
and using
\eqref{3.30}%
, we get%
\[
\gamma_{\alpha y}^{\prime\prime}(s^{\ast})=\zeta_{1}x(s)+\zeta_{2}%
\alpha(s)+\zeta_{3}y(s),
\]
where $\zeta_{1}=\frac{c\kappa^{\prime}}{b^{2}-2c^{2}}\left(  \frac{c-2\kappa
b}{2\kappa}\right)  ,\zeta_{2}=\frac{c\kappa^{\prime}}{b^{2}-2c^{2}}\left(
\frac{c-b-1}{2\kappa}\right)  ,\zeta_{3}=\frac{c\kappa^{\prime}}{b^{2}-2c^{2}%
}\left(  \frac{b\kappa-c}{2\kappa^{2}}\right)  .$%
\begin{equation}
y_{\alpha y}(s^{\ast})=-\gamma_{\alpha y}^{\prime\prime}-\frac{1}%
{2}\left\langle \gamma_{\alpha y}^{\prime\prime},\gamma_{\alpha y}%
^{\prime\prime}\right\rangle \gamma_{\alpha y}. \tag{3.32}\label{3.32}%
\end{equation}

By the help of equation
\eqref{3.32}%
, we obtain%
\begin{equation}
y_{\alpha y}(s^{\ast})=\omega_{1}x(s)+\omega_{2}\alpha(s)+\omega_{3}y(s),
\tag{3.33}\label{3.33}%
\end{equation}
where $\omega_{1}=-\zeta_{1},$ $\omega_{2}=-(\zeta_{2}+\frac{1}{2}\left(
2\zeta_{1}\zeta_{3}+\zeta_{2}^{2}\right)  ),$ $\omega_{3}=-(\zeta_{3}+\frac
{b}{2c}\left(  2\zeta_{1}\zeta_{3}+\zeta_{2}^{2}\right)  ).$

\textbf{ii)} The curvature $\kappa_{\gamma_{\alpha y}}(s^{\ast})$ of the
$\gamma_{\alpha y}(s^{\ast})$ is explicity obtained by%
\[
\kappa_{\gamma_{\alpha y}}(s^{\ast})=-\frac{c^{2}}{8\left(  b^{2}%
-2c^{2}\right)  ^{2}}\left(  \frac{\kappa^{\prime}}{\kappa^{{}}}\right)
^{2}\left(  \frac{\left(  c-2\kappa b\right)  \left(  b\kappa-c\right)
}{\kappa}+\left(  c-b-1\right)  ^{2}\right)  .
\]

\end{proof}

\begin{definition}
Let $x$ be unit speed spacelike curve lying on $\mathbf{Q}_{{}}^{2}$ with the
moving asymptotic orthonormal frame $\left\{  x,\alpha,y\right\}  .$ Then,
$x\alpha y-$smarandache curve of $x$ is defined by%
\begin{equation}
\gamma_{x\alpha y}(s^{\ast})=\frac{1}{\sqrt{2cc^{\ast}+b^{2}}}\left(
cx(s)+b\alpha(s)+c^{\ast}y(s)\right)  , \tag{3.34}\label{3.34}%
\end{equation}
where $c,c^{\ast},b\in%
\mathbb{R}
_{0}^{+}.$
\end{definition}

\begin{theorem}
Let $x$ be unit speed spacelike curve in $\mathbf{Q}_{{}}^{2}$ with the moving
asymptotic orthonormal frame $\left\{  x,\alpha,y\right\}  $ and cone
curvature $\kappa$ and let $\gamma_{x\alpha y}$ be $x\alpha y-$smarandache
curve with asymptotic orthonormal frame $\left\{  \gamma_{x\alpha y},\text{
}\alpha_{x\alpha y},\text{ }y_{x\alpha y}\right\}  .$ Then the following
relations hold:

i) The asymptotic orthonormal frame $\left\{  \gamma_{x\alpha y},\text{
}\alpha_{x\alpha y},\text{ }y_{x\alpha y}\right\}  $ of the $x\alpha
y-$smarandache curve $\gamma_{x\alpha y}$ is given as%
\begin{equation}%
\begin{bmatrix}
\gamma_{x\alpha y}\\
\alpha_{x\alpha y}\\
y_{x\alpha y}%
\end{bmatrix}
=%
\begin{bmatrix}
\frac{c}{\sqrt{2cc^{\ast}+b^{2}}} & \frac{b}{\sqrt{2cc^{\ast}+b^{2}}} &
\frac{c^{\ast}}{\sqrt{2cc^{\ast}+b^{2}}}\\
\rho_{1} & \rho_{2} & \rho_{3}\\
\sigma_{1} & \sigma_{2} & \sigma_{3}%
\end{bmatrix}%
\begin{bmatrix}
x\\
\alpha\\
y
\end{bmatrix}
\tag{3.35}\label{3.35}%
\end{equation}
where
\begin{align}
\eta &  =\sqrt{\left(  c-c^{\ast}\kappa(s)\right)  ^{2}-2b^{2}\kappa
(s)}\nonumber\\
\rho_{1}  &  =\frac{b\kappa(s)}{\eta},\rho_{2}=\frac{c-c^{\ast}\kappa(s)}%
{\eta},\rho_{3}=-\frac{b}{\eta}\tag{3.36}\label{3.36}\\
\xi_{1}  &  =\left(  \rho_{1}^{\prime}+\rho_{2}\kappa\right)  ,\xi_{2}%
=\rho_{2}^{\prime}+\rho_{1}+\rho_{3}\kappa,\xi_{3}=-\rho_{3}^{\prime}-\rho
_{2}\nonumber
\end{align}
and
\begin{align}
\sigma_{1}  &  =-\xi_{1}-\frac{c}{2\sqrt{2cc^{\ast}+b^{2}}}\left(  2\xi_{1}%
\xi_{3}+\xi_{2}^{2}\right)  ,\nonumber\\
\sigma_{2}  &  =-\xi_{2}-\frac{b}{2\sqrt{2cc^{\ast}+b^{2}}}\left(  2\xi_{1}%
\xi_{3}+\xi_{2}^{2}\right)  ,\tag{3.37}\label{3.37}\\
\sigma_{3}  &  =-\xi_{3}-\frac{c^{\ast}}{2\sqrt{2cc^{\ast}+b^{2}}}\left(
2\xi_{1}\xi_{3}+\xi_{2}^{2}\right)  .\nonumber
\end{align}

ii) The cone curvature $\kappa_{\gamma_{x\alpha y}}(s^{\ast})$ of the curve
$\gamma_{x\alpha y}$ is given by%
\begin{align}
\kappa_{\gamma_{yx\alpha}}(s^{\ast})  &  =\left(  b\left(  \frac{\kappa}{\eta
}\right)  ^{\prime}+\frac{c-c^{\ast}\kappa}{\eta}\kappa\right)  \left(
\left(  \frac{b}{\eta}\right)  ^{\prime}+\frac{c-c^{\ast}\kappa}{\eta}\right)
\nonumber\\
&  -\frac{1}{2}\left(  \frac{c-c^{\ast}\kappa}{\eta}\right)  ^{\prime}.
\tag{3.38}\label{3.38}%
\end{align}
where
\begin{equation}
s^{\ast}=\frac{1}{\sqrt{2cc^{\ast}+b^{2}}}\int\sqrt{\left(  c-c^{\ast}%
\kappa(s)\right)  ^{2}-2b^{2}\kappa(s)}ds,\text{ }b,c,c^{\ast}\in%
\mathbb{R}
_{0}^{+}. \tag{3.39}\label{3.39}%
\end{equation}

\end{theorem}

\begin{proof}
\textbf{i)} Differentiating the equation
\eqref{3.34}
with respect to $s$ and considering
\eqref{2.1}%
, we find%
\begin{equation}
\gamma_{x\alpha y}^{\prime}(s^{\ast})\frac{ds^{\ast}}{ds}=\frac{1}%
{\sqrt{2cc^{\ast}+b^{2}}}\left(  b\kappa\overrightarrow{x(s)}+\left(
c-c^{\ast}\kappa\right)  \overrightarrow{\alpha(s)}-\overrightarrow
{by(s)}\right)  .\tag{3.40}\label{3.40}%
\end{equation}
This can be written as follows%
\begin{equation}
\alpha_{x\alpha y}(s^{\ast})=\frac{b\kappa}{\eta}\overrightarrow{x(s)}%
+\frac{c-c^{\ast}\kappa}{\eta}\overrightarrow{\alpha(s)}-\frac{b}{\eta
}\overrightarrow{y(s)},\tag{3.41}\label{3.41}%
\end{equation}
or%
\begin{equation}
\alpha_{x\alpha y}(s^{\ast})=\rho_{1}\overrightarrow{x(s)}+\rho_{2}%
\overrightarrow{\alpha(s)}-\rho_{3}\overrightarrow{y(s)}.\tag{3.42}%
\label{3.42}%
\end{equation}
where
\begin{equation}
\frac{ds^{\ast}}{ds}=\frac{1}{\sqrt{2cc^{\ast}+b^{2}}}\sqrt{\left(  c-c^{\ast
}\kappa\right)  ^{2}-2b^{2}\kappa}\tag{3.43}\label{3.43}%
\end{equation}
Differentiating
\eqref{3.42}
and using
\eqref{3.43}%
, we get%
\[
\gamma_{xy\alpha}^{\prime\prime}(s^{\ast})=\xi_{1}x(s)+\xi_{2}\alpha
(s)+\xi_{3}y(s),
\]
where $\xi_{1}=\left(  \rho_{1}^{\prime}+\rho_{2}\kappa\right)  ,$ $\xi
_{2}=\rho_{2}^{\prime}+\rho_{1}+\rho_{3}\kappa,$ $\xi_{3}=-\rho_{3}^{\prime
}-\rho_{2}.$%
\begin{equation}
y_{x\alpha y}(s^{\ast})=-\gamma_{x\alpha y}^{\prime\prime}-\frac{1}%
{2}\left\langle \gamma_{x\alpha y}^{\prime\prime},\gamma_{x\alpha y}%
^{\prime\prime}\right\rangle \gamma_{x\alpha y}.\tag{3.44}\label{3.44}%
\end{equation}
By the help of equation
\eqref{3.44}%
, we obtain%
\begin{equation}
y_{xay}(s^{\ast})=\sigma_{1}x(s)+\sigma_{2}\alpha(s)+\sigma_{3}y(s),\tag{3.45}%
\label{3.45}%
\end{equation}
where $\sigma_{1}=-\xi_{1}-\frac{c}{2\sqrt{2cc^{\ast}+b^{2}}}\left(  2\xi
_{1}\xi_{3}+\xi_{2}^{2}\right)  ,$ $\sigma_{2}=-\xi_{2}-\frac{b}%
{2\sqrt{2cc^{\ast}+b^{2}}}\left(  2\xi_{1}\xi_{3}+\xi_{2}^{2}\right)  ,$
$\sigma_{3}=-\xi_{3}-\frac{c^{\ast}}{2\sqrt{2cc^{\ast}+b^{2}}}\left(  2\xi
_{1}\xi_{3}+\xi_{2}^{2}\right)  .$

\textbf{ii)} From $\kappa_{\gamma_{x\alpha y}}(s^{\ast})=-\frac{1}%
{2}\left\langle \gamma_{x\alpha y}^{\prime\prime},\gamma_{x\alpha y}%
^{\prime\prime}\right\rangle ,$ we have%
\[
\kappa_{\gamma_{yx\alpha}}(s^{\ast})=\left(  b\left(  \frac{\kappa}{\eta
}\right)  ^{\prime}+\frac{c-c^{\ast}\kappa}{\eta}\kappa\right)  \left(
\left(  \frac{b}{\eta}\right)  ^{\prime}+\frac{c-c^{\ast}\kappa}{\eta}\right)
-\frac{1}{2}\left(  \frac{c-c^{\ast}\kappa}{\eta}\right)  ^{\prime}.
\]

\end{proof}

\begin{theorem}
Let $x:I\rightarrow$ $\mathbf{Q}^{2}$\textbf{\ }$\subset E_{1}^{3}$ be a
spacelike curve in $\mathbf{Q}^{2}$ as follows%
\begin{equation}
x(s)=\frac{f_{s}^{-1}}{2}(f^{2}-1,2f,f^{2}+1), \tag{3.46}\label{3.46}%
\end{equation}
for some non constant function $f(s).$ Then we can write the following statements:

\begin{itemize}
\item If $x$ is a $x\alpha-$ smarandache curve, then the $x\alpha-$
smarandache curve $\gamma_{x\alpha}$ can be written as%
\begin{equation}
\gamma_{x\alpha}(s^{\ast})=(\frac{c}{b}-f_{s}^{-1}f_{ss}^{{}})x(s)+\left(
f,1,f\right)  . \tag{3.47}%
\end{equation}

\item If $x$ is a $xy-$ smarandache curve, then the $xy-$ smarandache curve
$\gamma_{xy}$ can be written as%
\begin{equation}
\gamma_{xy}(s^{\ast})=\frac{1}{\sqrt{2bc}}\left(  (c-\frac{1}{2}f_{s}%
^{-2}f_{ss}^{2})x(s)+f_{s}^{-1}f_{ss}^{{}}\left(  f,1,f\right)  -f_{s}\left(
1.0,1\right)  \right)  . \tag{3.48}%
\end{equation}

\item If $x$ is a $\alpha y-$ smarandache curve, then the $\alpha y-$
smarandache curve $\gamma_{\alpha y}$ can be written as%
\begin{align}
\gamma_{\alpha y}(s^{\ast})  &  =\left(  -f_{s}^{-1}f_{ss}^{{}}-\frac{b}%
{2c}f_{s}^{-2}f_{ss}^{2}\right)  x(s)+\left(  1+\frac{b}{c}f_{s}^{-1}%
f_{ss}^{{}}\right)  \left(  f,1,f\right) \tag{3.49}\\
&  -\frac{b}{c}f_{s}\left(  1.0,1\right)  .\nonumber
\end{align}

\item If $x$ is a $x\alpha y-$ smarandache curve, then the $x\alpha y-$
smarandache curve $\gamma_{x\alpha y}$ can be written as%
\[
\gamma_{x\alpha y}(s^{\ast})=\frac{1}{\sqrt{2cc^{\ast}+b^{2}}}\left(
\begin{array}
[c]{c}%
\left(  c-bf_{s}^{-1}f_{ss}^{{}}-\frac{c^{\ast}}{2}f_{s}^{-2}f_{ss}%
^{2}\right)  x(s)\\
+\left(  b+c^{\ast}f_{s}^{-1}f_{ss}^{{}}\right)  \left(  f,1,f\right)
-f_{s}\left(  1.0,1\right)
\end{array}
\right)  ,
\]

\end{itemize}

where $c,b,c^{\ast}\in%
\mathbb{R}
_{0}^{+}.$
\end{theorem}

\begin{proof}
It is obvious from
\eqref{3.1}%
,
\eqref{3.14}%
,
\eqref{3.23}%
,
\eqref{3.34}
and
\eqref{3.46}%
.
\end{proof}

We can give the following example to hold special Smarandache curves in the
null cone $\mathbf{Q}^{2}.$ $x\alpha$, $xy$, $\alpha y,$ and $x\alpha y$
special smarandache curves of $x$ curves are given in Figure 1 A, C, E, G, I,
respectively. These figures rotated in three dimensions are also given in
Figure 1 B, D, F, H, J, respectively.

\begin{example}
The curve%
\[
x(s)=\left(  \frac{\cosh s}{2}-\frac{1}{\cosh s},\tanh s,\frac{\cosh s}%
{2}\right)
\]
is spacelike in $\mathbf{Q}^{2}$ with arc length parameter $s$. Also, the
shape of the $x$ curve is given as follows
\end{example}

\begin{figure}[ptb]
\centering
\subfloat[The curve $x $ \label{fig1a}] {\
\includegraphics[width=6cm,height=3cm]{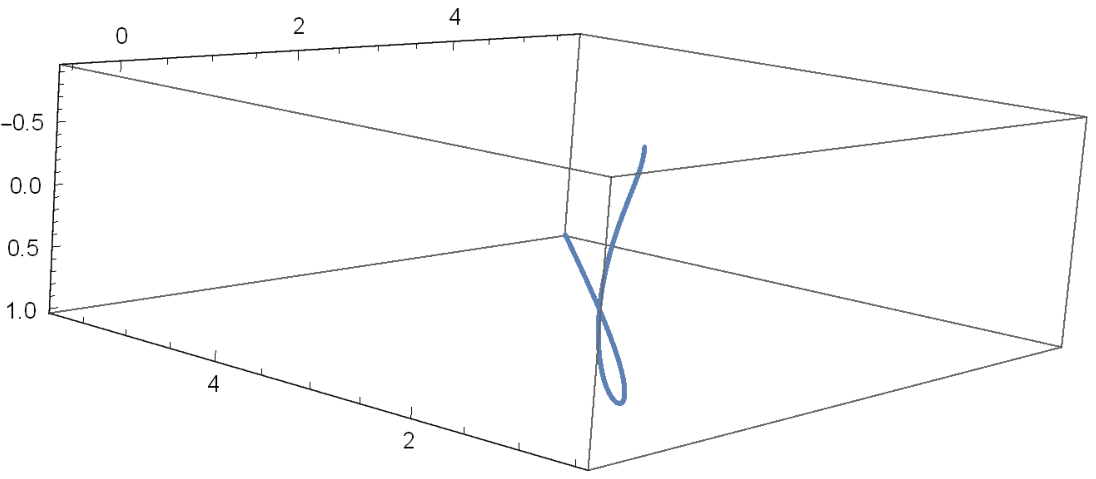} } \hspace*{.7cm}
\subfloat[The rotated surface of  curve $x$ \label{fig1b}] {\
\includegraphics[width=6cm,height=3cm]{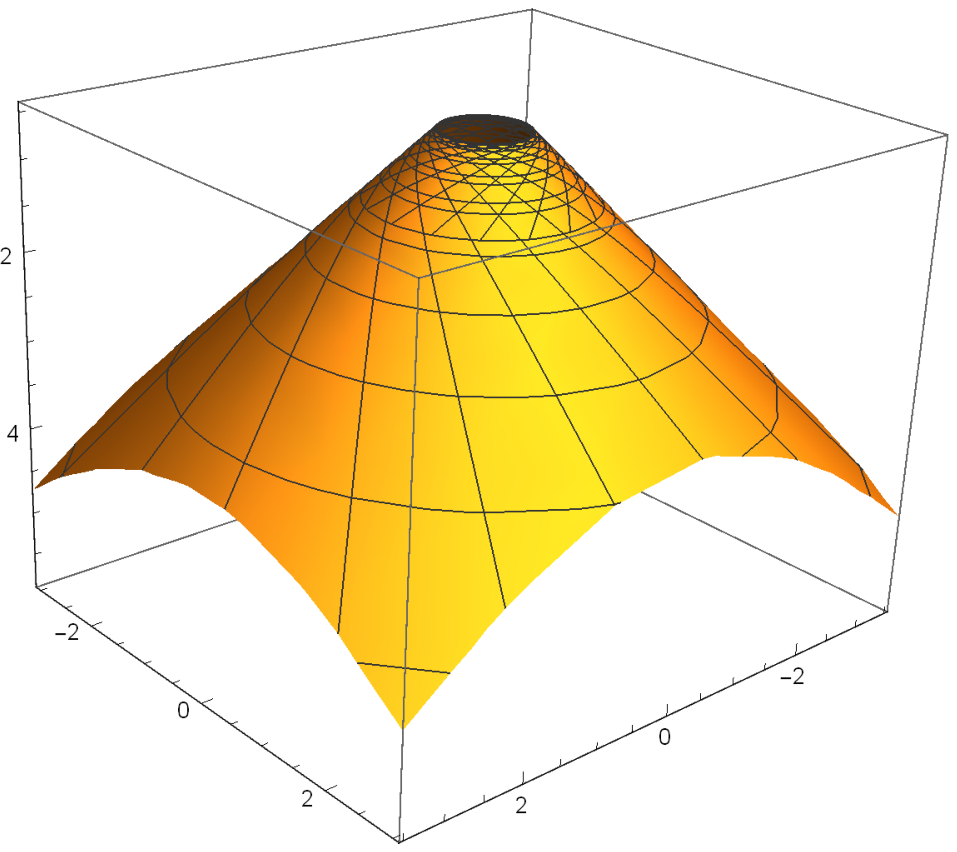} } \newline
\par
\subfloat[$ x\alpha- $smarandache curve \label{fig2a}] {\
\includegraphics[width=6cm,height=3cm]{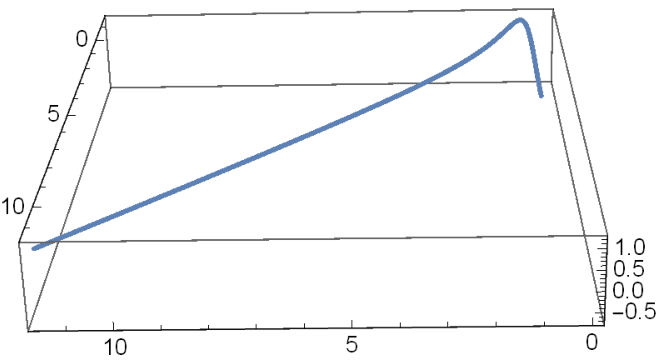} } \hspace*{.7cm}
\subfloat[$ x\alpha $- smarandache surface. \label{fig2b}] {\
\includegraphics[width=6cm,height=3cm]{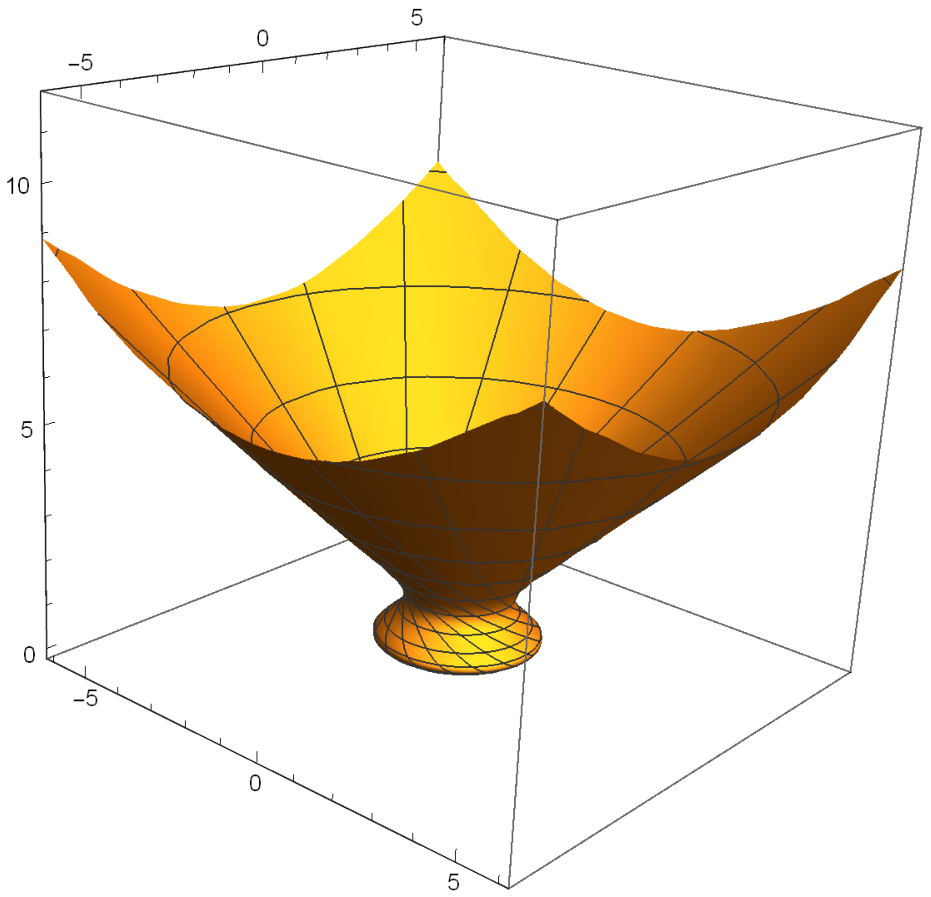} } \newline
\par
\subfloat[$ xy- $smarandache curve \label{fig3a}] {\
\includegraphics[width=6cm,height=3cm]{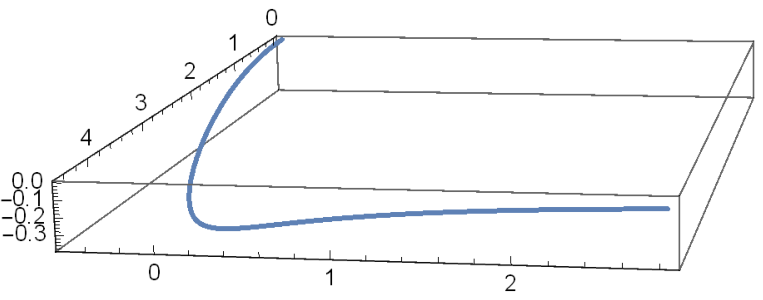} } \hspace*{.7cm}
\subfloat[$ xy- $smarandache surface. \label{Fig3b}] {\
\includegraphics[width=6cm,height=3cm]{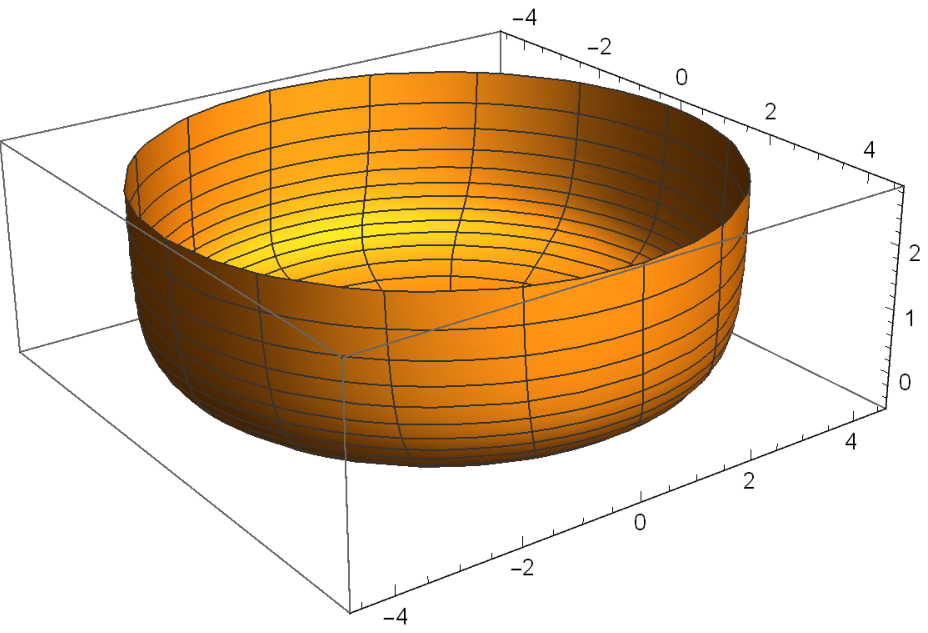} } \newline
\par
\subfloat[$\alpha$y$- $smarandache curve.\label{fig4a}] {\
\includegraphics[width=6cm,height=3cm]{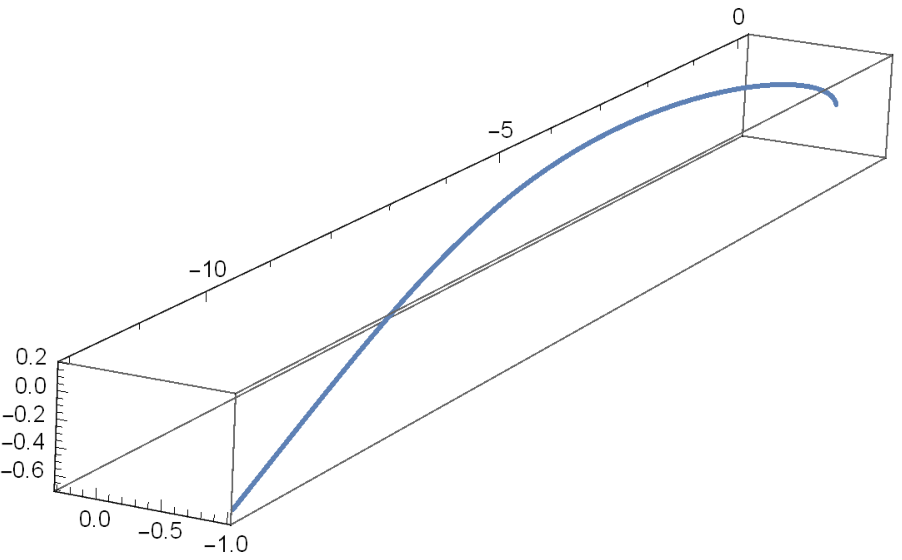} } \hspace*{.7cm}
\subfloat[$\alpha$y$- $smarandache surface. \label{fig4b}] {\
\includegraphics[width=6cm,height=3cm]{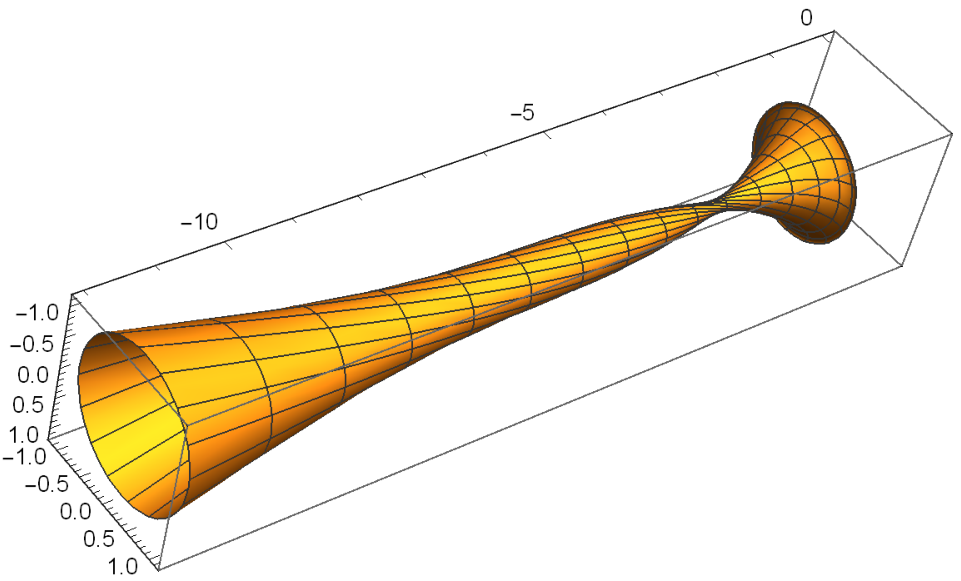} }
\par
\subfloat[$x\alpha$y$- $smarandache curve.\label{fig5a}] {\
\includegraphics[width=6cm,height=3cm]{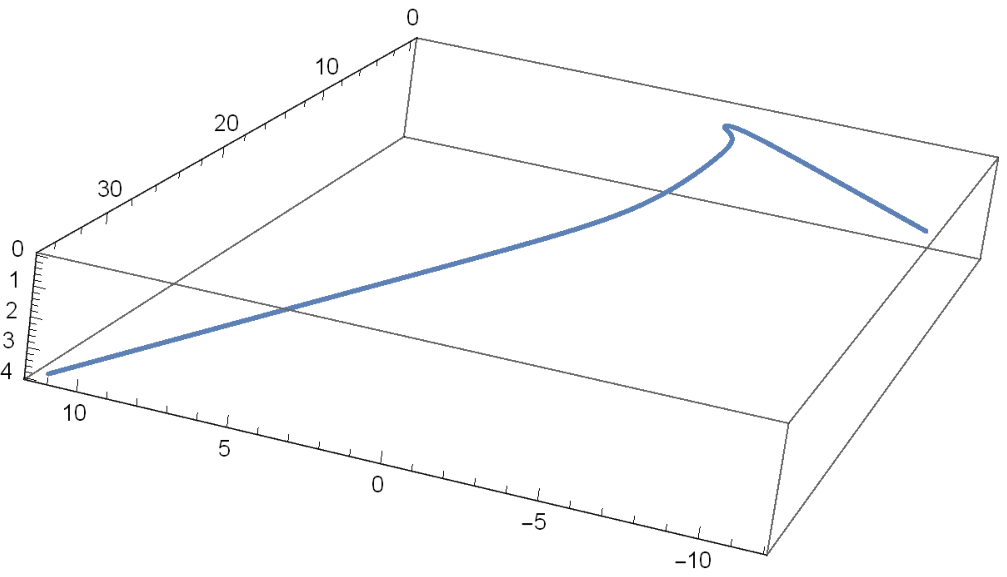} } \hspace*{.7cm}
\subfloat[$ x\alpha$y$- $smarandache surface. \label{fig5b}] {\
\includegraphics[width=6cm,height=3cm]{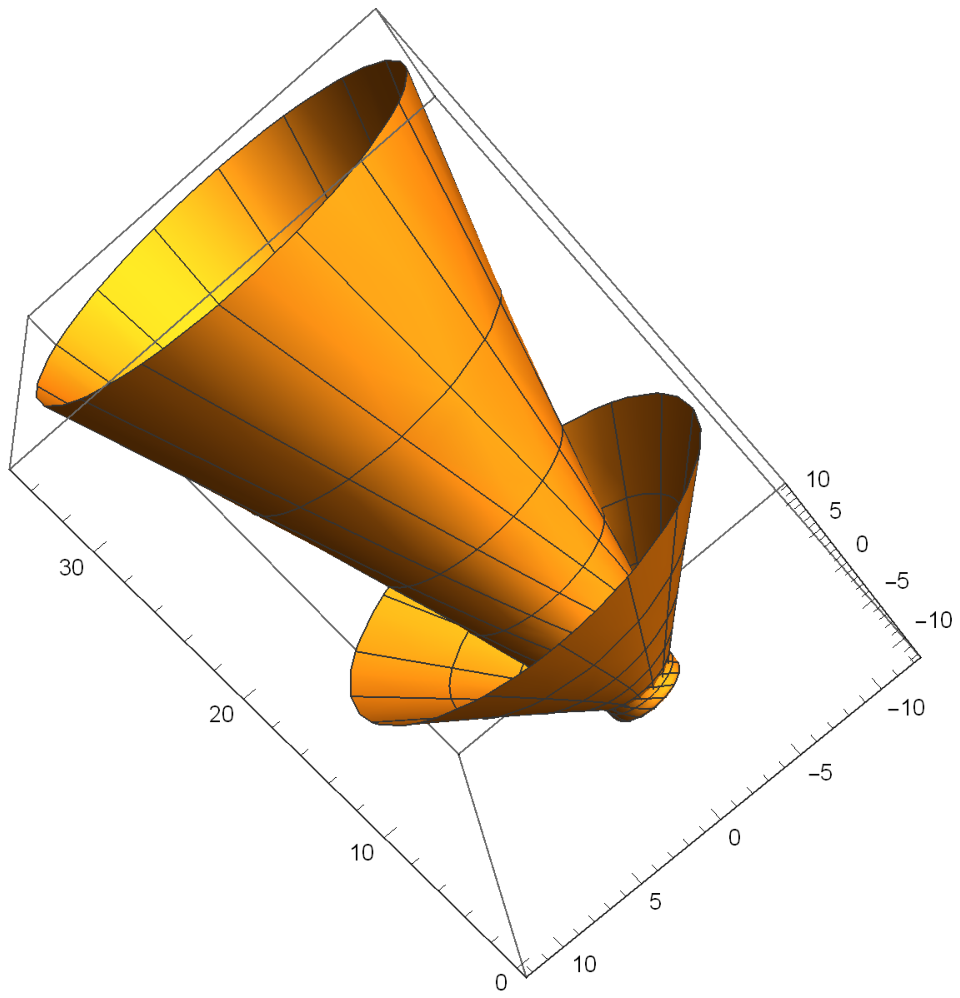} }\caption{}\label{Fig1}%
\end{figure}


Then we can write the smarandache curves of the $x$-curve as follows:

\textbf{i)} $x\alpha-$ smarandache curve $\gamma_{x\alpha}$ is given by%
\[
\gamma_{x\alpha}(s)=\left(
\begin{array}
[c]{c}%
d\left(  \frac{\cosh s}{2}-\frac{1}{\cosh s}\right)  +\frac{\sinh s}{2}%
-\frac{\tanh s}{\cosh s},\\
d\tanh s+\frac{1}{\cosh^{2}s},\\
d\frac{\cosh s}{2}+\frac{\sinh s}{2}%
\end{array}
\right)
\]

\textbf{ii)} \textbf{\ }$xy-$ smarandache curve $\gamma_{xy}$ is given by%
\[
\gamma_{xy}(s)=\left(
\begin{array}
[c]{c}%
m\cosh s-n\tanh s\sinh s-d\left(  \frac{1c+\tanh^{2}s}{\cosh s}\right)  ,\\
c\tanh s\left(  e-\frac{\tanh^{2}s}{2}\right)  ,\\
d\left(  \left(  c-\frac{\tanh^{2}s}{2}\right)  \frac{\cosh s}{2}-\frac
{1}{\cosh s}\right)
\end{array}
\right)
\]

\textbf{iii)} $\alpha y-$ smarandache curve $\gamma_{\alpha y}$ is given by%
\[
\gamma_{\alpha y}(s)=\left(
\begin{array}
[c]{c}%
\left(  \frac{\cosh s}{2}-\frac{1}{\cosh s}\right)  \left(  1-e\tanh s\right)
\tanh s+\sinh s\left(  1+d\tanh s\right)  \\
-d\cosh s,\\
-\tanh^{2}s\left(  1+d\tanh s\right)  +d\tanh s,\\
\frac{\sinh s}{2}-d\cosh s-e\sinh s\tanh s
\end{array}
\right)
\]

iv) $x\alpha y-$ smarandache curve $\gamma_{x\alpha y}$ is given by%
\[
\gamma_{x\alpha y}(s)=\left(
\begin{array}
[c]{c}%
m\sinh s-n\cosh s+d\sinh s\tanh s-\frac{c}{\cosh s}-\frac{b\sinh s-c^{\ast
}\tanh^{2}s}{\cosh^{2}s},\\
m+n\tanh s+k\tanh^{2}s+l\tanh^{3}s,\\
m\cosh s+d\sinh s+c\tanh s\sinh s,
\end{array}
\right)
\]

where $m,n,d,c,k,l,e\in%
\mathbb{R}
_{0}^{+}.$

\end{document}